\documentclass[11pt]{article}

\usepackage[T1]{fontenc}
\usepackage[utf8]{inputenc}
\usepackage[letterpaper,margin=1in]{geometry}
\usepackage{amsmath,amssymb}
\usepackage{graphicx}
\usepackage{booktabs}
\usepackage{float}
\usepackage[numbers,sort&compress]{natbib}
\usepackage[colorlinks,allcolors=blue]{hyperref}
\usepackage[expansion=false]{microtype}

\title{Verified parameter estimation for differential-equation models:\\
NIST reference problems and the theophylline flip-flop}

\author{Hamidreza Semiyari\thanks{Luminary Quantum Institute LLC, Silver Spring,
Maryland, USA. \texttt{MIE@mitrainstituteofeducation.org}. ORCID
\href{https://orcid.org/0000-0001-5531-737X}{0000-0001-5531-737X}.}}

\date{\today}

\begin{document}
\maketitle

\begin{abstract}
Parameters fitted to differential-equation models carry two kinds of error:
statistical error from noisy data, and numerical error from the solver and the
optimizer that compute the fit. Standard fitting tools report the first and leave
the second unexamined, and they return one optimum without saying whether another
parameter set fits the data equally well. We describe SIRPY, a differential-equation
solver that checks every solution it returns by substituting it back into the
equation, and we test it in two settings. First, on eleven nonlinear-regression
reference problems from the US National Institute of Standards and Technology
(NIST), each written as a first-order differential equation and fitted by least
squares with every model curve supplied by SIRPY, the estimates agree with NIST's
certified values to 7.1 to 8.8 significant digits on ten problems, with all 6{,}281
solves verified; the eleventh problem exposed a limitation of the residual check,
which we report and explain. Second, on the classic theophylline pharmacokinetic
data, SIRPY found two parameter sets whose predicted curves agree to within
$1.7\times10^{-11}$~mg/L and matched the reference fit from the R package
\texttt{nlme} to every printed digit; the data alone cannot decide which rate is
absorption and which is elimination, and physiology must. We state precisely what
the verification certifies, what it does not, and what remains to be done to turn a
residual check into a proven bound on estimation error.
\end{abstract}


\section{Introduction}

A parameter estimated by fitting a differential-equation model to data is the output
of three computations: a solver that produces the model curve for each trial
parameter set, an optimizer that moves the parameters toward the best fit, and a
statistical calculation that turns the curvature of the fit into standard errors and
intervals. The second and third are documented in every statistics textbook. The
first is usually treated as exact. It is not: a numerical solver returns an
approximation, the optimizer stops at a tolerance, and neither error is reported
alongside the estimate. In pharmacokinetics, where a fitted clearance sets a
maintenance dose, this silent error sits inside a number that is used as a fact.

A second, older problem is that nonlinear fits can have more than one optimum. The
one-compartment oral model is the standard example: the absorption and elimination
rate constants can be exchanged without changing the predicted concentration curve
\citep{wagner1993,yanez2011}, so the data alone cannot tell which is which. Fitting
tools return whichever optimum the starting point leads to.

This paper reports a solver built to make the first problem visible and a test of it
against the second. SIRPY is a solver for ordinary differential equations based on
the power-series (Parker--Sochacki) method \citep{parker1996,carothers2004}, extended
by the author to boundary-value and Volterra integral problems
\citep{semiyari2015,semiyari2016}. Its distinguishing feature is that every solution
it returns is checked by substitution into the equation; the result carries a
verdict, and a solution that fails the check is not used. We ask two questions. On
problems with answers certified to eleven digits by the US National Institute of
Standards and Technology \citep{nist_strd}, how close does a least-squares fit come
when every model curve is a verified SIRPY solution? And on the theophylline data
distributed with the R package \texttt{nlme} \citep{pinheiro2000,boeckmann1994}, does
the method find both optima and say what the data cannot decide?

We also report what the first test taught us about the verification itself. A
residual check is a necessary condition, not a proof of accuracy, and one NIST
problem shows why.

\section{Methods}

\subsection{What SIRPY verifies}

For an initial-value problem $y' = f(x, y; \theta)$, $y(x_0) = y_0$, SIRPY returns a
piecewise polynomial approximation $\tilde y(x)$ on $[x_0, x_1]$ together with a
verdict. The verdict is VERIFIED when the residual
$\tilde y'(x) - f(x, \tilde y(x); \theta)$, evaluated on a dense set of points across
the interval, stays below a stated tolerance and the initial condition is met.
Otherwise the verdict is NOT VERIFIED with the location and size of the violation, or
a named special point (a blow-up or branch point) when the solution ceases to exist.
The check is a statement about the returned function: it satisfies the equation to
the stated tolerance. It is not, by itself, a bound on $|\tilde y - y|$, which
requires a stability constant for the equation; we return to this in
Section~\ref{sec:eckerle} and the Discussion.

The solver's construction is the author's published extension of the
Parker--Sochacki method \citep{semiyari2015}; its implementation is proprietary to
Luminary Quantum Institute LLC. Everything reported here depends only on the solver's
output contract: the returned function, its evaluation at requested points, and the
verdict.

\subsection{NIST reference problems}\label{sec:nist-methods}

The NIST Statistical Reference Datasets include 27 nonlinear least-squares problems
with parameter values and residual sums of squares certified to 11 significant
digits, each with two starting points, the first deliberately far from the solution
\citep{nist_strd}. We selected the eleven whose models are solutions of a first-order
ordinary differential equation in the predictor $x$ (Table~\ref{tab:models}),
spanning NIST's lower, average and higher difficulty levels.

\begin{table}[H]
\centering
\small
\begin{tabular}{llcll}
\toprule
Problem & Level & $n$ & Model & Differential equation used \\
\midrule
Misra1a  & lower   & 14  & $b_1(1-e^{-b_2 x})$                & $y' = b_2(b_1 - y)$, $y(0)=0$ \\
Misra1b  & lower   & 14  & $b_1(1-(1+b_2 x/2)^{-2})$          & $y' = b_1 b_2 (1-y/b_1)^{3/2}$, $y(0)=0$ \\
Misra1c  & average & 14  & $b_1(1-(1+2b_2 x)^{-1/2})$         & $y' = b_1 b_2 (1-y/b_1)^{3}$, $y(0)=0$ \\
Misra1d  & average & 14  & $b_1 b_2 x/(1+b_2 x)$              & $y' = b_2(b_1-y)^2/b_1$, $y(0)=0$ \\
Chwirut2 & lower   & 54  & $e^{-b_1 x}/(b_2+b_3 x)$           & $y' = -(b_1 + b_3/(b_2+b_3 x))\,y$ \\
DanWood  & lower   & 6   & $b_1 x^{b_2}$                      & $y' = b_2 y/x$ \\
BoxBOD   & higher  & 6   & $b_1(1-e^{-b_2 x})$                & $y' = b_2(b_1 - y)$, $y(0)=0$ \\
Rat42    & higher  & 9   & $b_1/(1+e^{b_2-b_3 x})$            & $y' = b_3 y(1-y/b_1)$ \\
MGH10    & higher  & 16  & $b_1 e^{b_2/(x+b_3)}$              & $y' = -b_2 y/(x+b_3)^2$ \\
Eckerle4 & higher  & 35  & $(b_1/b_2)\,e^{-\frac{1}{2}((x-b_3)/b_2)^2}$ & $y' = -(x-b_3)\,y/b_2^2$ \\
Bennett5 & higher  & 154 & $b_1(b_2+x)^{-1/b_3}$              & $y' = -y/(b_3(b_2+x))$ \\
\bottomrule
\end{tabular}
\caption{The eleven NIST problems and the differential equations solved for them.
Where the initial condition is not $y(0)=0$, the solve starts at the smallest data
$x$ with the initial value taken from the model formula, and the curve beyond it is
SIRPY's solution.}
\label{tab:models}
\end{table}

For each problem the fit minimized the sum of squared residuals with the
Levenberg--Marquardt implementation in SciPy \citep{scipy2020}, started from both
NIST starting points, with tolerances of $10^{-14}$ and at most 600 function
evaluations per pass. Every function evaluation was one SIRPY solve over the data
range followed by evaluation of the returned function at the data points; a solve
whose verdict was not VERIFIED, or that exceeded five seconds, was discarded and the
evaluation treated as infeasible. As an independent cross-check, the same
least-squares fit was run with the closed formula and no differential equation.
Agreement with NIST is reported as significant digits,
$-\log_{10}|\hat b - b^\star|/|b^\star|$, for the worst parameter of each problem.
Every run reported in this paper, in this section and the next, used SIRPY 0.67.0 on
a single Apple computer; the scripts, the raw results files and the NIST inputs are
provided (Data and code availability).

\subsection{Theophylline}\label{sec:theo-methods}

The theophylline data of Boeckmann, Sheiner and Beal \citep{boeckmann1994}, as
distributed in R, give plasma concentrations of twelve subjects after an oral dose.
We used subject 1 (dose 4.02~mg/kg, eleven samples over 24.4~h, including the
pre-dose sample) and the one-compartment model with first-order absorption,
\begin{equation}
G' = -k_a G,\qquad C' = \frac{k_a}{V}\,G - k_e C,\qquad G(0)=D,\; C(0)=0,
\end{equation}
with volume $V = CL/k_e$ and Gaussian measurement error; the error variance is
concentrated out of the likelihood, so maximum likelihood coincides with least
squares. Parameters were estimated on the log scale within the box
$k_a, k_e \in [0.01, 10]$~h$^{-1}$, $CL \in [0.005, 0.5]$~L/h/kg.

Three stages were run. (i) A scan of the $(k_a, k_e)$ plane on a $29\times29$
logarithmic grid, with $CL$ optimized at each node, to map the likelihood surface;
841 solves. (ii) From the best grid node on each side of the line $k_a = k_e$, a
certified maximum-likelihood fit: Newton iterations on the log-likelihood with every
function value a verified SIRPY solve, followed by a check that the score vanishes
and the observed information is negative definite, and 95\% profile-likelihood
intervals for each parameter. (iii) Cross-checks: the Bateman closed form evaluated
at the eleven sample times against the SIRPY solution, and an independent
least-squares fit of the closed form with SciPy. The reference values are those
printed by \texttt{example(SSfol)} in the \texttt{nlme} package
\citep{pinheiro2000}, which fits the same subject with the same model:
$\log k_e = -2.9196$, $\log k_a = 0.5752$, $\log CL = -3.9159$, residual standard
error 0.732 on 8 degrees of freedom.

\section{Results}

\subsection{NIST problems}\label{sec:nist-results}

\begin{figure}[H]
\centering
\includegraphics[width=0.9\textwidth]{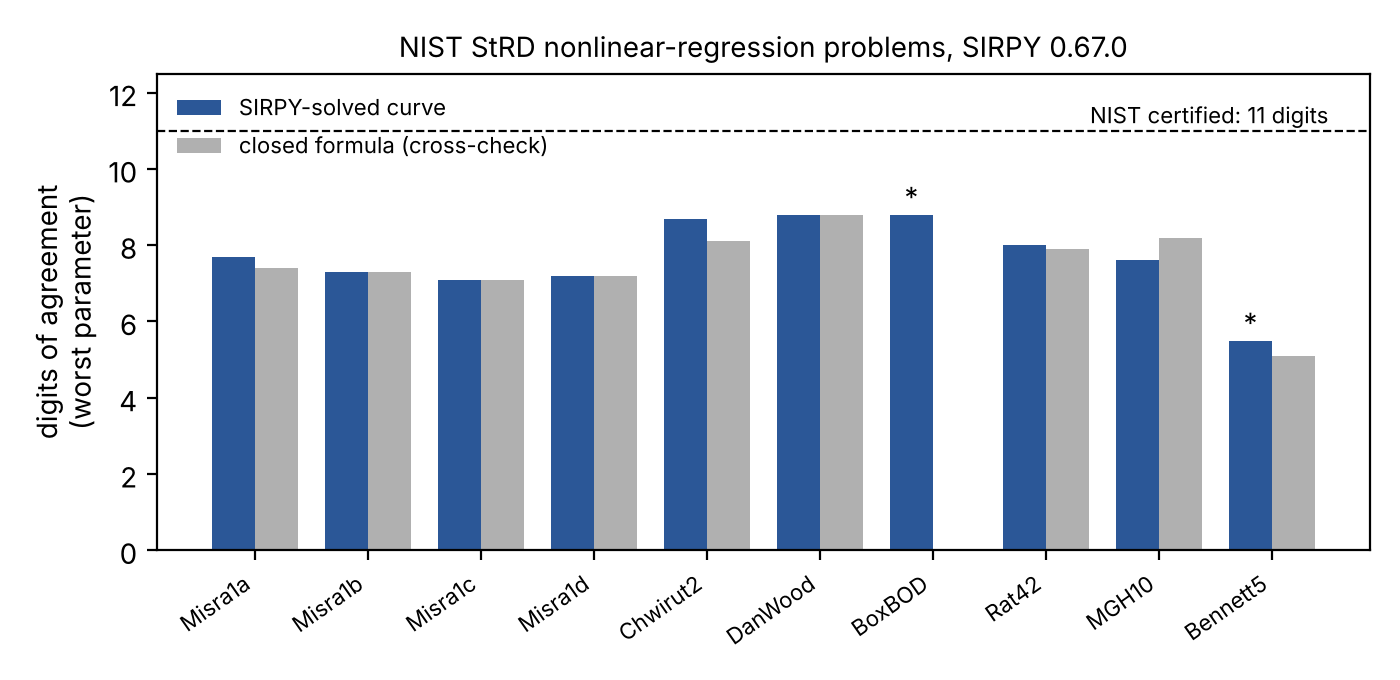}
\caption{Significant digits of agreement with NIST's certified values, worst
parameter, for the ten problems with a sound differential-equation formulation. Grey
bars are the closed-formula cross-check; a star marks a problem where the two NIST
starting points did not reach the same optimum. Eckerle4 is discussed separately.}
\label{fig:nist}
\end{figure}

Table~\ref{tab:nist} and Figure~\ref{fig:nist} summarize the eleven problems. On ten
of them the SIRPY-driven fit agrees with NIST's certified parameters to between 7.1
and 8.8 significant digits on the worst parameter, and the residual sum of squares
agrees to 8.2 to 11.7 digits. Every solve used in those ten fits passed verification
(6{,}281 solves in total; three further solves on BoxBOD exceeded the time limit and
were discarded). The closed-formula cross-check reached the same answers to within
one digit, except on BoxBOD, where the closed-formula fit started from NIST's first
point ($b_1 = 1$, $b_2 = 1$) converged to a wrong optimum ($b_1 = 172.5$,
$b_2 = 110.9$, a nearly flat curve) while the SIRPY-driven fit found the certified
one from its second start. Bennett5, which NIST lists as its most ill-conditioned
problem, reached 5.5 digits by both routes; the limit there is the optimizer's
stopping precision, not the solver.

\begin{table}[H]
\centering
\footnotesize
\setlength{\tabcolsep}{4pt}
\begin{tabular}{llcccccccc}
\toprule
 &  &  & Digits & Digits & Digits & Same optimum & & & \\
Problem & Level & $n$ & (SIRPY) & (RSS) & (closed form) & from both starts & Solves & Verified & Seconds \\
\midrule
Misra1a  & lower   & 14  & 7.7    & 10.4 & 7.4    & yes & 115   & 115   & 0.9 \\
Misra1b  & lower   & 14  & 7.3    & 11.1 & 7.3    & yes & 133   & 133   & 4.9 \\
Misra1c  & average & 14  & 7.1    & 8.2  & 7.1    & yes & 88    & 88    & 0.9 \\
Misra1d  & average & 14  & 7.2    & 10.4 & 7.2    & yes & 69    & 69    & 0.5 \\
Chwirut2 & lower   & 54  & 8.7    & 11.2 & 8.1    & yes & 99    & 99    & 2.5 \\
DanWood  & lower   & 6   & 8.8    & 11.7 & 8.8    & yes & 69    & 69    & 0.5 \\
BoxBOD   & higher  & 6   & 8.8    & 10.4 & $-2.3$ & no  & 125   & 124   & 57.6 \\
Rat42    & higher  & 9   & 8.0    & 11.0 & 7.9    & yes & 101   & 101   & 1.2 \\
MGH10    & higher  & 16  & 7.6    & 9.5  & 8.2    & yes & 1{,}535 & 1{,}535 & 25.2 \\
Eckerle4 & higher  & 35  & $-6.3$ & 11.4 & 10.3   & no  & 2{,}053 & 2{,}053 & 81.1 \\
Bennett5 & higher  & 154 & 5.5    & 11.5 & 5.1    & no  & 3{,}947 & 3{,}947 & 49.7 \\
\bottomrule
\end{tabular}
\caption{Agreement with NIST's certified values. ``Digits'' is the number of
significant digits shared with the certified value, for the worst parameter; NIST
certifies 11. ``Verified'' counts the solves that passed SIRPY's check and were used.
SIRPY 0.67.0, one run, one computer.}
\label{tab:nist}
\end{table}

\subsection{Eckerle4: what a residual check does not prove}\label{sec:eckerle}

On Eckerle4 the fitted residual sum of squares matched NIST to 11.4 digits and two
parameters to 7 and 10 digits, while $b_1$ came out as $3.1\times10^{6}$ instead of
1.554. The cause is the formulation in Table~\ref{tab:models}, not the data. The
equation $y' = -(x-b_3)y/b_2^2$ is linear and homogeneous, so it does not contain the
amplitude $b_1$ at all; $b_1$ enters only through the initial value at $x_0 = 400$,
where the model is about $1.2\times10^{-35}$. Starting from a value that small and
growing by a factor of $e^{79}$ to the peak, the solver returned a curve about
$3\times10^{6}$ times too small at the peak, and reported it VERIFIED, because the
residual check is absolute: at a scale of $10^{-35}$ to $10^{-10}$ every residual is
tiny in absolute terms. The optimizer then compensated by inflating $b_1$.

Two lessons follow, and both are being acted on. The verification tolerance must be
relative to the scale of the solution, or the certificate must state the scale at
which the check was made; and a solve whose solution spans tens of orders of
magnitude should be refused rather than certified. For the benchmark itself the
remedy is to fit $u = \ln y$, whose equation $u' = -(x-b_3)/b_2^2$ has a well-scaled
initial value; we report the problem as it was run rather than after the fact.

\subsection{Theophylline}\label{sec:theo-results}

\begin{figure}[H]
\centering
\includegraphics[width=0.8\textwidth]{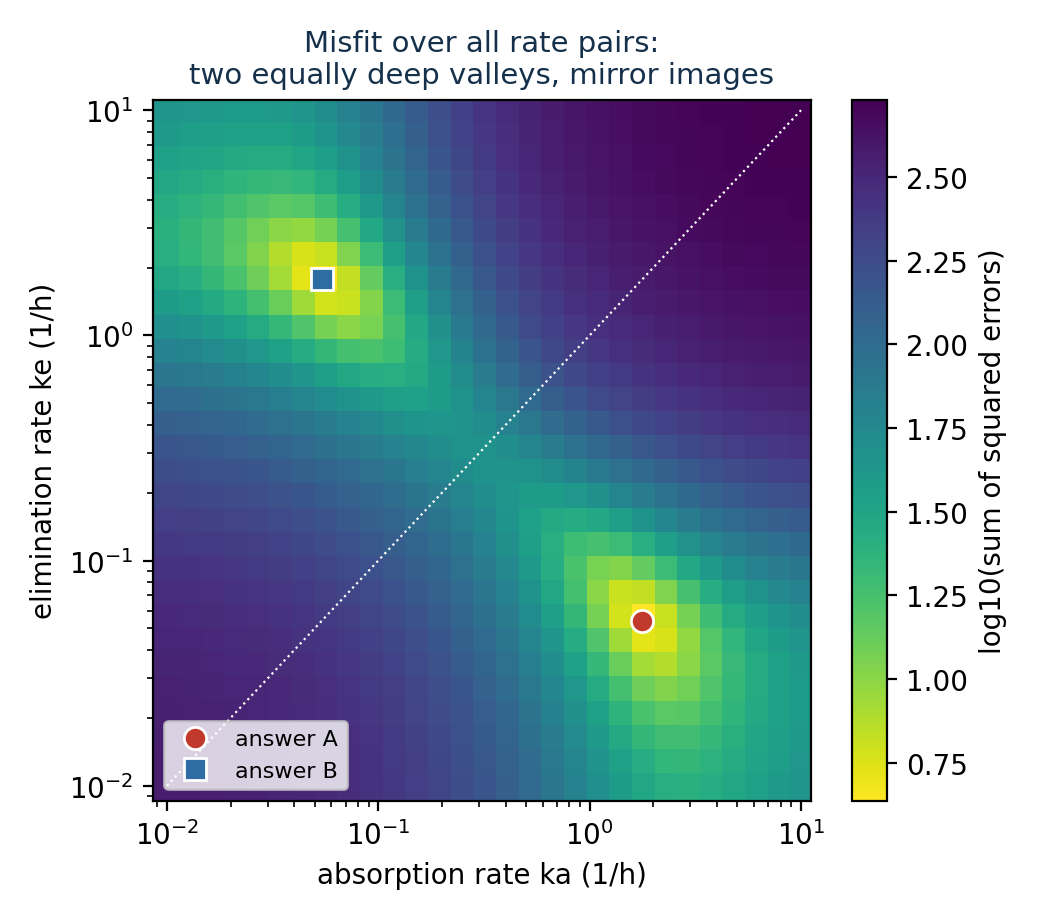}
\caption{The misfit surface over the $(k_a, k_e)$ plane from the 841-solve scan, with
clearance optimized at each node. The two basins are mirror images across
$k_a = k_e$.}
\label{fig:landscape}
\end{figure}

The scan of the $(k_a, k_e)$ plane (Figure~\ref{fig:landscape}) shows two basins
separated by the line $k_a = k_e$. Started from the best grid node in each basin, the
certified estimator converged and certified both answers (Table~\ref{tab:theo}): at
each reported point the largest score component was $2.8\times10^{-10}$, the observed
information was negative definite in every direction, and all three parameters were
identified. Within each basin, nine starting points spread across that half of the
box reached the same single mode, which is evidence of uniqueness there rather than a
proof of it.

\begin{table}[H]
\centering
\small
\begin{tabular}{lll}
\toprule
 & Answer A & Answer B \\
\midrule
$k_a$ (h$^{-1}$)   & 1.777 [1.351, 2.362]      & 0.05395 [0.04835, 0.06094] \\
$k_e$ (h$^{-1}$)   & 0.05395 [0.04835, 0.06094] & 1.777 [1.351, 2.362] \\
$CL$ (L/h/kg)      & 0.01992 [0.01810, 0.02186] & 0.01992 [0.01810, 0.02186] \\
$V = CL/k_e$ (L/kg) & 0.369                     & 0.0112 \\
Elimination half-life & 12.8 h                  & 0.39 h \\
Sum of squared residuals & 4.28600902           & 4.28600902 \\
SIRPY solves (all verified) & 18{,}438          & 18{,}305 \\
Max difference from closed form at sample times (mg/L) & $9.1\times10^{-13}$ & $9.2\times10^{-13}$ \\
\bottomrule
\end{tabular}
\caption{Two certified maximum-likelihood fits to subject 1. Brackets are 95\%
profile-likelihood intervals. Answer A agrees with the \texttt{nlme} reference fit to
every printed digit ($\log k_e = -2.9196$, $\log k_a = 0.5752$,
$\log CL = -3.9159$).}
\label{tab:theo}
\end{table}

\begin{figure}[H]
\centering
\includegraphics[width=0.8\textwidth]{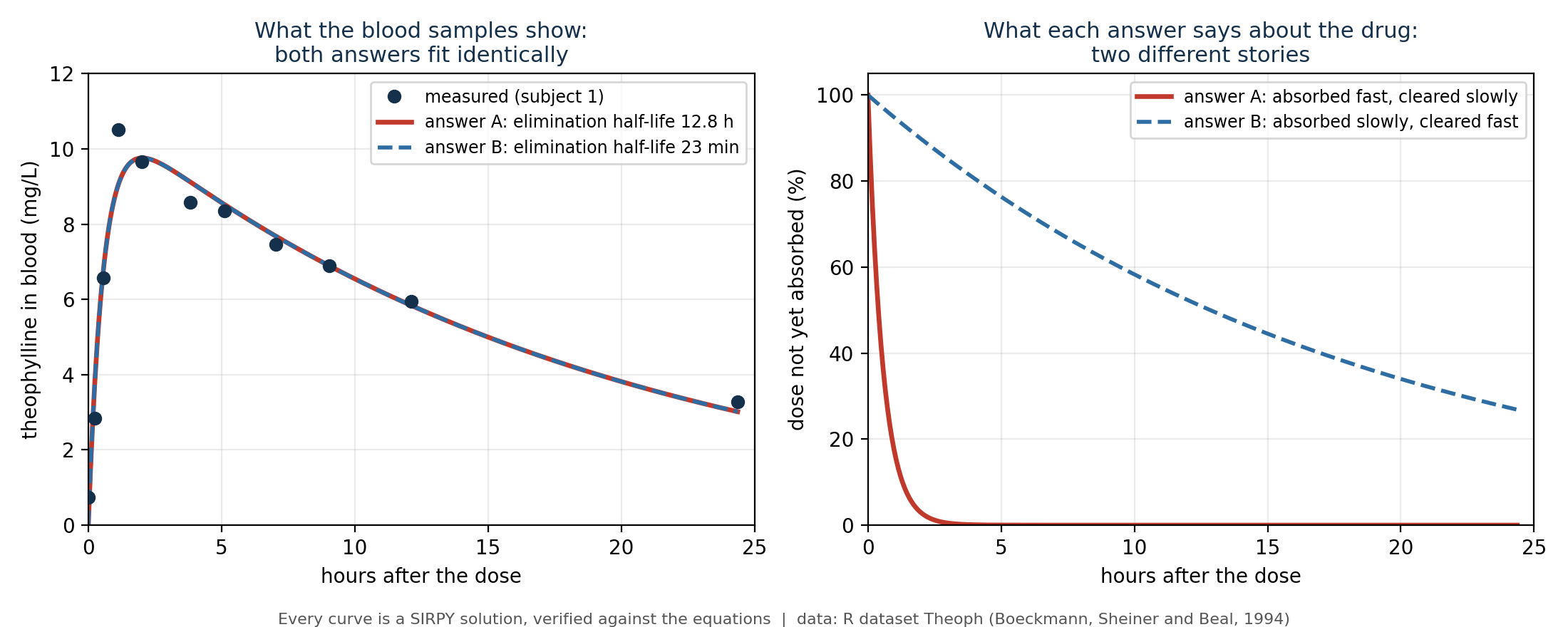}
\caption{Subject 1's measured concentrations with the two certified fits, which
coincide to $1.7\times10^{-11}$~mg/L.}
\label{fig:flipflop}
\end{figure}

The two predicted concentration curves differ by at most $1.7\times10^{-11}$~mg/L
over the sampled range (Figure~\ref{fig:flipflop}). Clearance, and therefore the area
under the curve, is the same in both; what differs is which rate is absorption.
Answer B requires a distribution volume of 0.011~L/kg, smaller than the plasma volume
of roughly 0.04~L/kg, so physiology selects A. The data cannot. The independent
closed-form least-squares fit agrees with each certified answer to a relative
$2\times10^{-8}$. In all, 37{,}584 SIRPY solves were made in this analysis (841 in
the scan, 18{,}438 and 18{,}305 in the two fits) and every one passed verification.

\section{Discussion}

Two things are certified in these results and one is not. Certified: every model
curve used in every fit satisfied its differential equation to the stated tolerance,
and each theophylline optimum is a point where the score vanishes and the observed
information is negative definite, with profile intervals computed rather than
approximated. Also established, by an outside reference rather than by the solver: on
ten NIST problems the fits agree with certified values to 7 to 9 digits, which is at
or above what published fitting software reports on the same problems
\citep{mccullough1998,mccullough1999}. Not certified: a bound on the numerical error
inside each estimate. The Eckerle4 case shows that a residual check can pass a
solution that is wrong by six orders of magnitude when the solution's scale makes the
absolute residual meaningless, and the theophylline certified fits carry no statement
of how far solver error could have moved the optimum. Turning the residual into a
bound on solution error, carrying that bound through the likelihood derivatives into
the estimates and intervals, and testing for a unique optimum within a computed
region, are the next steps; they are the subject of the author's current research
program.

The cost structure is also informative. A single SIRPY solve on these problems takes
10 to 50~ms, and a least-squares fit needs about a hundred; the certified
theophylline fits with profile intervals needed about 18{,}000 solves each, because
the estimator's derivatives are taken by finite differences. Exact derivatives from
the sensitivity equations would cut that by an order of magnitude and are also what
an error bound requires.

Finally, the theophylline result is a reminder that identifiability is a property of
the data and the model, not of the software. A tool that reports one optimum hides
the choice between A and B; a tool that reports both, with the physiological argument
alongside, lets the modeler make it. That matters whenever the fitted parameters are
used to predict a situation the data did not cover, such as an intravenous dose or a
formulation with a different absorption rate.

\section*{Data and code availability}

The NIST problem files are public \citep{nist_strd}. The theophylline data are
distributed with R (\texttt{Theoph}). The benchmark script, the raw results file
\texttt{nist\_results.json}, the theophylline scan and fit outputs, and the figures
are available at \url{https://github.com/sirpy26/SIRPY}. SIRPY is proprietary
software of Luminary Quantum Institute LLC; its solver output and verdicts are
reproducible from the released version stated in each results file.

\section*{Acknowledgments and disclosure}

The author developed SIRPY and owns it through Luminary Quantum Institute LLC. AI
tools were used for drafting the analysis scripts and editing the text; all
computations were run by the author on his own computer and all numbers were checked
by him.

\bibliography{references}

\end{document}